\documentclass[11pt,a4paper]{article}

\usepackage{amsmath}
\usepackage{amsopn}
\usepackage{amsfonts}
\usepackage{amsthm}
\usepackage{fullpage}

\numberwithin{equation}{section}

\theoremstyle{definition}

\theoremstyle{remark}

\newtheoremstyle{mytheorem}{0.5cm}{0.2cm}{\slshape}{ }{\bfseries}{.}{ }{}
\theoremstyle{mytheorem}
\newtheorem{theorem}{Theorem}[section]

\newtheorem{lemma}[theorem]{Lemma}
\newtheorem{corollary}[theorem]{Corollary}

\renewcommand{\P}{\mathbf{P}}

\newcommand{\standardspace}[1]{\mathbb{#1}}      
\newcommand{\R}{\standardspace{R}}

\newcommand{\E}{{\mathbf{E}}}

\newcommand{\thf}{\frac{1}{2}}

\DeclareMathOperator{\RV}{RV}

\usepackage{epsfig,psfrag}

\newcommand{\D}{\mathcal{D}}

\begin{document}

\title{On the domain of attraction for the lower tail in Wicksell's
  corpuscle problem}

\author{S. K\"otzer \quad I. Molchanov \\\\
\normalsize{Department of Mathematical Statistics and Actuarial
Sciences,}\\
  \normalsize{University of Bern, Switzerland}}
\date{}

\maketitle

\begin{abstract}

We consider the classical Wicksell corpuscle problem with
spherical particles in $\mathbb{R}^n$ and investigate the shapes
of lower tails of distributions of `sphere radii' in
$\mathbb{R}^n$ and `sphere radii' in a $k$-dimensional section
plane $E_k$. We show in which way the domains of attraction are
related to each other.

\end{abstract}

\section{Wicksell's corpuscle problem}
\label{sec:wicksell}

Suppose that a collection of spheres is randomly scattered in the
Euclidean space $\R^n$ and we intersect the collection with a
$k$-dimensional section plane $E_k$, $k\in\{1,\ldots,n-1\}$. Then,
under suitable model assumptions, the relation between the radii
distribution $F$ of spheres in $\R^n$ and the radii distribution
$F^{(n,k)}$ of corresponding spheres in the section plane is given
by Wicksell's integral equation,
\begin{align}
 \label{wicksellrelation}
  F^{(n,k)}(x)=1-\frac{(n-k)}{M_{n-k}}\int_x^{\infty}u(u^2-x^2)^
  {\frac{n-k-2}{2}}(1-F(u))\,du\,,\quad x>0\,.
\end{align}
where $M_{n-k}\in(0,\infty)$ denotes the $(n-k)$th moment of $F$.
Relation (\ref{wicksellrelation}) has been discovered in the
mid-twenties by S.D. Wicksell \cite{wick25} for the special case
$n=3$ and $k=2$. A proof of (\ref{wicksellrelation}) can be found
in \cite[Ch. 4.3]{schnweil00}. Note that the density of
$F^{(n,k)}$ always exists and is given by
\begin{align}
 \label{wicksellrelationdensity}
  f^{(n,k)}(x)=\frac{x\,(n-k)}{M_{n-k}}\int_x^\infty
  (u^2-x^2)^{\frac{n-k-2}{2}}\,dF(u)\,,\quad x>0\,,
\end{align}
which can easily be seen by differentiating
(\ref{wicksellrelation}). From (\ref{wicksellrelation}) and
(\ref{wicksellrelationdensity}) we see that the representations of
$F^{(n,k)}$ and $f^{(n,k)}$ only depend on the difference $n-k$.
In the following we will therefore write $F^{(r)}$ for
$F^{(n,k)}$, $f^{(r)}$ for $f^{(n,k)}$, and $M_r$ for $M_{n-k}$,
where $r=n-k$.

In recent years there has been some effort in creating a theory of
\textit{stereology of extremes}. Drees and Reiss \cite{dreereis92}
have investigated the shapes of upper tails of distributions of
`sphere radii' and `circle radii' that are connected by
(\ref{wicksellrelation}). In this paper we investigate the lower
tail behaviour.

\section{Extreme value theory and regular variation}
\label{sec:extreme-value-theory-regular-variation}

Let $\{X_i, \,i\geq1\}$ be a sequence of independent and
identically distributed real valued random variables with common
distribution function $F$. Denote the sample minimum by
\begin{displaymath}
  W_n=\min(X_1,\ldots,X_n)\,,\quad n\geq 1\,.
\end{displaymath}
If there are sequences of normalizing constants $a_n>0$ and
$b_n\in\R$ such that the normalized minima $(W_n-b_n)/a_n$
converge in distribution to a nondegenerate distribution function
$H$, then $F$ lies in the domain of attraction of $H$. We denote
this by $F\in\D(H)$. There are only three types of possible
limiting distributions $H$, see e.g. \cite{resn87}, namely
\begin{displaymath}
  H_{i,\alpha}(x)=
  \begin{cases}
    1-\exp(-(-x)^{-\alpha}) & \quad x<0\,, \quad i=1\,, \\
    1-\exp(-x^\alpha) & \quad x>0\,, \quad i=2\,,
  \end{cases} \quad \alpha>0\,,
\end{displaymath}
and
\begin{displaymath}
  H_3(x)=1-\exp(-e^x)\,,\quad x\in\R\,.
\end{displaymath}
$H_{1,\alpha}$, $H_{2,\alpha}$ and $H_3$ are called Fr\'echet,
Weibull and Gumbel distributions respectively. Recall the
following condition for $F\in\D(H_{2,\alpha})$. Let
$\eta=\inf\{x\in\R:\,F(x)>0\}$ denote the lower endpoint of $F$.
Then $F\in\D(H_{2,\alpha})$ if and only if $\eta>-\infty$ and for
all $x>0$,
\begin{eqnarray}
  \label{min-nec-suff-condition2}
  \lim_{s\downarrow0}\frac{F(\eta+xs)}{F(\eta+s)}=x^{\alpha}\,.
\end{eqnarray}

In the following we recall some basic facts from the theory of
regular variation, see e.g. \cite{sene76}. A function $R$ is
regularly varying at infinity with exponent $\rho\in\R$ if it is
real-valued, positive and measurable on $[x_0,\infty)$, for some
$x_0>0$, and if for each $x>0$
\begin{displaymath}
 \lim_{t\to\infty}\frac{R(tx)}{R(t)}=x^\rho\,.
\end{displaymath}
We then write $R\in \RV_\infty(\rho)$. A function $R(\cdot)$ is
regularly varying at 0 if $R(\frac{1}{x})$ is regularly varying at
infinity with exponent $(-\rho)$ and this is denoted by $R\in
\RV_0(\rho)$. If $\rho=0$, $R$ is said to be slowly varying. A
function $R(\cdot)$ is regularly varying if and only if it can be
written in the form
\begin{align}
 \label{product-form-of-regvarying-functions}
 R(x)=x^\rho L(x)\,,
\end{align}
where $\rho\in\R$ and $L(\cdot)$ is slowly varying. Note that the
sum of two slowly varying functions is again slowly varying. For
the investigation of the lower tail in Wicksell's corpuscle
problem we need the following lemmas.

\begin{lemma}
 \label{lemma:sum-of-slowly-varying-functions}
  If $R_1\in \RV_0(\beta_1)$ and $R_2\in \RV_0(\beta_2)$ then $R_1+R_2\in
  \RV_0(\min\{\beta_1,\beta_2\})$.
\end{lemma}

\begin{lemma}[cf. Theorem 2.7 in \cite{sene76}]
 \label{lemma: integral slow variation}
 Let $L$ be slowly varying at $0$ on $(0,\infty)$ and bounded on
 each finite subinterval of $(0,\infty)$. Suppose that for $\beta>0$
 the integral $\int_{\beta}^\infty t^{\delta}f(t)\,dt$
 is well-defined for some given real function $f$ and a given
 number $\delta\geq 0$. Then as $u\downarrow 0$
 \begin{align*}
  \int_{\beta}^\infty f(t)L(ut)\,dt\sim L(u)\int_{\beta}^\infty f(t)\,dt
 \end{align*}
 for $\delta>0$, and also for $\delta=0$ provided that $L$ is
 non-increasing on $(0,\infty)$.
\end{lemma}

\begin{lemma}[cf. Exercise 1.1.1. in \cite{resn87}]
 \label{property: Gumbel}
 Let $F\in\D(H_3)$ and $\eta=\inf\{x\in\R:\,F(x)>0\}>-\infty$.
 Then for all $n\geq 1$,
 \begin{displaymath}
  \lim_{x\downarrow\eta}\,(x-\eta)^{-n}F(x)=0\,.
 \end{displaymath}
\end{lemma}

\begin{lemma}
 \label{lemma}
  Let $f(\cdot,\cdot)$ and $g(\cdot,\cdot)$ be positive functions
  such that
  \begin{displaymath}
   \int_0^\omega f(s,t)\,dt<\infty\,,\quad\int_0^\omega g(s,t)\,dt<\infty
  \end{displaymath}
  for some $\omega\in(0,\infty]$. Furthermore assume for $s\leq
  t<\omega$
  \begin{displaymath}
   \lim_{s\uparrow\omega}\frac{f(s,t)}{g(s,t)}=c\,,\quad
   c\in[0,\infty]\,.
  \end{displaymath}
  Then
  \begin{displaymath}
   \lim_{s\uparrow\omega}\frac{\int_s^\omega f(s,t)\,dt}{\int_s^\omega
   g(s,t)\,dt}=c\,.
  \end{displaymath}
\end{lemma}

\section{Tail behaviour in Wicksell's corpuscle problem}
\label{sec:tail-behav-wicks}

In recent years the main theoretical tools for estimating
stereologically the tail of a particle size distribution have been
stability properties of the domain of attraction. Consider the
spherical Wicksell corpuscle problem with corresponding
distribution functions $F$ and $F^{(r)}$, $r\geq 1$, as described
in Section \ref{sec:wicksell}.

Note that there is no need to investigate the behaviour for
$F\in\D(H_{1,\alpha})$ since in this case
$\eta=\inf\{x\in\R:\,F(x)>0\}=-\infty$. Let us first assume that
$r=1$ and $F\in\D(H_{2,\alpha})$. The radius of the spheres in
$\R^n$ is denoted by $\xi$. We have

\begin{theorem}
 \label{minwik-behaviour}
  Let $F\in\D(H_{2,\alpha})$ and $\eta=0$. Then $F^{(1)}\in\D(H_{2,\beta})$,
  where
  \begin{displaymath}
   \beta=\begin{cases}
         2, & \alpha>1\,,\\
         \alpha+1, & 0<\alpha\leq1\,.
         \end{cases}
  \end{displaymath}
  If $\eta>0$, then
  $F^{(1)}\in\D(H_{2,2})$.
\end{theorem}

\begin{proof}
  Let us first assume that $\eta=0$. It is easily seen from
  (\ref{wicksellrelationdensity}) that
  \begin{align}
   \label{proof:lowertail1}
   F^{(1)}(t) & =\frac{1}{M_1}\left[\int_0^t u\,dF(u)+\int_t^\infty
   (u-\sqrt{u^2-t^2})\,dF(u)\right]
   =:\frac{1}{M_1}[I_1(t)+I_2(t)]\,.
  \end{align}
  Using integration by parts in the first summand of
  (\ref{proof:lowertail1}) we get
  \begin{displaymath}
   I_1(t)=tF(t)-\int_0^t F(u)\,du\,.
  \end{displaymath}
  For all $\alpha>0$ and $x>0$ we show that
  \begin{align}
   \label{proof:lowertail2}
    \frac{I_1(tx)}{I_1(t)}\rightarrow x^{\alpha+1}
    \quad\quad \textrm{as}\,\,\,t\downarrow0\,,
  \end{align}
  i.e. $I_1\in \RV_0(\alpha+1)$. For that consider
  \begin{align*}
   \frac{I_1(tx)}{I_1(t)} & =\frac{txF(tx)-\int_0^{tx} F(u)\,du}
   {tF(t)-\int_0^t F(u)\,du}=\frac{xF(tx)}{F(t)}\cdot
   \frac{1-(txF(tx))^{-1}\int_0^{tx} F(u)\,du}
   {1-(tF(t))^{-1}\int_0^{t} F(u)\,du}\,.
  \end{align*}
  Define $G(s)=s^{-2}F(\frac{1}{s})\in \RV_\infty(-\alpha-2)$ and
  apply Karamata's theorem \cite[Th.~0.6]{resn87} to
  obtain
  \begin{align*}
   \lim_{t\downarrow0}\frac{\int_0^t F(u)\,du}{tF(t)}=
   \lim_{z\to\infty}\frac{\int_z^\infty
   G(s)\,ds}{zG(z)}=\frac{1}{\alpha+1}\,,
  \end{align*}
  which proves (\ref{proof:lowertail2}). \\
  For the second summand $I_2(t)$ we consider three cases. \\
  Case I, $\alpha>1$: Since $\alpha>1$, $\E(\xi^{-1})<\infty$.
  Rewriting $I_2(t)$ as
  \begin{displaymath}
   I_2(t)=\E\left((\xi-\sqrt{\xi^2-t^2})I_{\{\xi>t\}}\right)=t^2\E\left((\xi+\sqrt{\xi^2-t^2})^{-1}I_{\{\xi>t\}}\right)
  \end{displaymath}
  and applying the dominated convergence theorem we get
  \begin{align}
   \label{proof:lowertail4}
   \lim_{t\downarrow0}\E\left((\xi+\sqrt{\xi^2-t^2})^{-1}
   I_{\{\xi>t\}}\right)=\thf\E(\xi^{-1})\in(0,\infty)\,,
  \end{align}
  which implies $I_2\in \RV_0(2)$. By Lemma \ref{lemma:sum-of-slowly-varying-functions},
  $F^{(1)}\in \RV_0(2)$, whence $F^{(1)}\in\D(H_{2,2})$. \\
  Case II, $\alpha\in(0,1)$: We show that $I_2\in
  \RV_0(\alpha+1)$, which implies $F^{(1)}\in \RV_0(\alpha+1)$ by
  Lemma \ref{lemma:sum-of-slowly-varying-functions}. For that we
  use again integration by parts to arrive at
  \begin{align}
   \label{proof:lowertail3}
    I_2(t)=-tF(t)-\int_t^\infty
    F(u)(1-\frac{u}{\sqrt{u^2-t^2}})\,du\,.
  \end{align}
  Writing $F(t)=t^{\alpha}L(t)$ where $L$ is slowly varying at
  $0$, and substituting $u=tv$ in the integral of
  (\ref{proof:lowertail3}) yields that
  \begin{displaymath}
   I_2(t)=-t^{\alpha+1}L(t)+t^{\alpha+1}\int_1^\infty
   L(tv)v^{\alpha}g(v)\,dv\,,
  \end{displaymath}
  where
  \begin{align}
   \label{proof:lowertail6}
    g(v)=\frac{v}{\sqrt{v^2-1}}-1\,,\quad v>1\,,
  \end{align}
  is a probability density function. For $\delta\in(0,1-\alpha)$ we have $\int_1^\infty
  t^{\delta}t^{\alpha}g(t)\,dt<\infty$ and therefore Lemma
  \ref{lemma: integral slow variation} implies
  \begin{displaymath}
   h(t):=\frac{\int_1^\infty L(tv)v^{\alpha}g(v)\,dv}{c L(t)}\rightarrow1\,
   \quad\quad \textrm{as}\,\,\,t\downarrow0\,,
  \end{displaymath}
  where $c:=\int_1^\infty v^{\alpha}g(v)\,dv\in (1,\infty)$. Thus
  \begin{displaymath}
   I_2(t)=t^{\alpha+1}L(t)[c h(t)-1]\in \RV_0(\alpha+1)\,.
  \end{displaymath}
  Case III, $\alpha=1$: We show that $I_2\in \RV_0(2)$. Choose $A>1$ such that $c(A):=\int_1^A
  vg(v)\,dv>1$. From (\ref{proof:lowertail3}) we get
  \begin{align}
   \label{proof:lowertail5}
    I_2(t)=t\left[\int_1^A F(vt)g(v)\,dv-F(t)+\int_A^\infty F(vt)g(v)\,dv\right]\,,
  \end{align}
  where $g$ is given by (\ref{proof:lowertail6}). Consider
  \begin{displaymath}
   J_1(t)=\int_1^A F(vt)g(v)\,dv-F(t)=t\int_1^A L(vt)v
   g(v)\,dv-tL(t)\,,
  \end{displaymath}
  where $L$ is slowly varying at $0$. Using the uniform convergence theorem
  for slowly varying functions we get
  \begin{displaymath}
   h^\ast(t):=\frac{\int_1^A L(tv)v g(v)\,dv}{c(A) L(t)}\rightarrow1\,
   \quad\quad \textrm{as}\,\,\,t\downarrow0\,,
  \end{displaymath}
  whence
  \begin{displaymath}
   J_1(t)=tL(t)[c(A) h^\ast(t)-1]\in \RV_0(1)\,.
  \end{displaymath}
  The second summand in (\ref{proof:lowertail5}) can be written as
  \begin{displaymath}
   J_2(t)=\int_A^\infty F(vt)g(v)\,dv=\int_A^\infty \P(v^{-1}\xi\leq
   t)g(v)\,dv=\P(\zeta^{-1}\xi\leq t, \zeta>A)\,,
  \end{displaymath}
  where $\zeta$ is a random variable with probability density function $g$
  independent of $\xi$. Furthermore
  \begin{displaymath}
   \P(\zeta^{-1}\xi\leq t, \zeta>A)=
   \P(\zeta^{-1}\xi\leq t \mid \zeta>A)\P(\zeta>A)=
   F_{\tilde{\zeta}^{-1}\xi}(t)\P(\zeta>A)\,,
  \end{displaymath}
  where $\tilde{\zeta}$ is a random variable with probability density function
  $\tilde{g}: (A,\infty)\rightarrow (0,\infty)$
  given by $\tilde{g}(u)=c^{-1}g(u)$ for $c=\int_A^\infty
  g(u)\,du$. It is easy to check that
  $F_{\tilde{\zeta}^{-1}}\in \RV_0(1)$. By assumption
  $F_{\xi}\in \RV_0(1)$ and since $\xi$ and $\tilde{\zeta}^{-1}$
  are independent we have $F_{\tilde{\zeta}^{-1}\xi}\in
  \RV_0(1)$ (see e.g. \cite{embrgold80}, Theorem 3 and Corollary therein)
  and therefore $J_2\in \RV_0(1)$. Hence
  \begin{displaymath}
   I_2(t)=t[J_1(t)+J_2(t)]\in \RV_0(2)\,.
  \end{displaymath}
  Now let $\eta>0$. For sufficiently small t we get from
  (\ref{proof:lowertail1})
  \begin{displaymath}
    F^{(1)}(t)  =\frac{1}{M_1}\int_\eta^\infty
    (u-\sqrt{u^2-t^2})\,dF(u)=\frac{t^2}{M_1}\E(\xi+\sqrt{\xi^2-t^2})^{-1}\,.
  \end{displaymath}
  For all $\alpha>0$ we have $\E(\xi^{-1})<\infty$, whence
  $F^{(1)}\in\D(H_{2,2})$, using the same argument as in
  (\ref{proof:lowertail4}). The proof is complete.
\end{proof}

The case of a Gumbel limiting distribution is considered in the
next theorem.

\begin{theorem}
\label{minwik-behaviour-Gumbel}
  If $F\in\D(H_3)$, then $F^{(1)}\in\D(H_{2,2})$.
\end{theorem}

\begin{proof}
 The case $\eta>0$ is covered by Theorem
 \ref{minwik-behaviour}. So let $\eta=0$.
 From (\ref{proof:lowertail1}) we see that
 \begin{displaymath}
  F^{(1)}(t)=\frac{t}{M_1}\left[F_{\xi/\eta_1}(t)-F_{\xi/\eta_2}(t)\right]\,,
  \end{displaymath}
  where $\eta_1$ is a random variable with density $g$ given by
  (\ref{proof:lowertail6}) and $\eta_2$ is uniformly distributed on $[0,1]$
  and $\eta_1,\eta_2,\xi$ are independent. It suffices to show that
  $F_{\xi/\eta_1}\in \RV_0(1)$ and
  $F_{\xi/\eta_2}(t)/F_{\xi/\eta_1}(t)\rightarrow 0$ as
  $t\downarrow0$.\\
  Let us first prove that $F_{\xi/\eta_1}\in \RV_0(1)$. For that
  we write
  \begin{align}
   \label{proof:Gumbel}
   F_{\xi/\eta_1}(t) & =\P(\xi\leq t\eta_1,\eta_1\leq 1/t)
   +\P(\xi\leq t\eta_1,\eta_1> 1/t) \nonumber \\
   & =\frac{1}{t}\int_t^1 F(u)g(u/t)\,du+\frac{1}{t}\int_1^\infty
   F(u)g(u/t)\,du \nonumber \\
   & =:J^\star_1(t)+J^\star_2(t)\,.
  \end{align}
  For the second summand in (\ref{proof:Gumbel}) we can apply
  Lemma \ref{lemma} and get
  \begin{displaymath}
   \frac{J^\star_2(tx)}{J^\star_2(t)}=x^{-1}\frac{\int_{1/t}^\infty
   F(tw)(\frac{w}{\sqrt{w^2-x^2}}-1)\,dw}{\int_{1/t}^\infty
   F(tw)(\frac{w}{\sqrt{w^2-1}}-1)\,dw}\longrightarrow x^{-1}x^2=x\quad\quad \textrm{as}\,\,\,
   t\downarrow0\,.
  \end{displaymath}
  Since $g\in \RV_\infty(-2)$ the first summand in
  (\ref{proof:Gumbel}) can be written as
  \begin{displaymath}
   J^\star_1(t)=t\int_0^1 F(u)u^{-2}L(u/t)\,du\,,
  \end{displaymath}
  where $L(1/t)$ is slowly varying at 0.\\
  Because of Lemma \ref{property: Gumbel}, $\int_0^1
  u^{-\delta-2}F(u)\,du<\infty$ for some $\delta>0$ and therefore
  we can apply \cite[Theorem 2.7]{sene76} and get
  \begin{displaymath}
   h^{\ast\ast}(t)=\frac{\int_0^1
   F(u)u^{-2}L(u/t)\,du}{L(1/t)\int_0^1 F(u)u^{-2}\,du}\rightarrow
   1\quad\quad \textrm{as}\,\,\,t\downarrow0\,,
  \end{displaymath}
  i.e.
  \begin{displaymath}
   J^\star_1(t)=tL(1/t)h^{\ast\ast}(t)\int_0^1 F(u)u^{-2}\,du\in
   \RV_0(1)\,.
  \end{displaymath}
  Next we show that $F_{\xi/\eta_2}(t)/F_{\xi/\eta_1}(t)\rightarrow 0$ as
  $t\downarrow0$. Since $t^{-2}F_{\xi/\eta_1}(t)\in \RV_0(-1)$, we
  have
  \begin{displaymath}
   \lim_{t\downarrow0} t^{-2}F_{\xi/\eta_1}(t)=\infty\,.
  \end{displaymath}
  Furthermore Lemma \ref{property: Gumbel} implies
  \begin{displaymath}
   \lim_{t\downarrow0} t^{-2}F_{\xi/\eta_2}(t)=\lim_{t\downarrow0}
   t^{-2}\int_0^1 F(tu)\,du=0\,.
  \end{displaymath}
  The result follows.
\end{proof}

The sectioning of the system of spheres with a $k$-dimensional
plane can be obtained by an iterated intersection procedure. First
we intersect the system with an $(n-1)$-dimensional plane
$E_{n-1}$, then the obtained system of spheres in $E_{n-1}$ is
intersected with an $(n-2)$-dimensional plane and so on. This
leads to

\begin{corollary}
 \label{minwik-behaviour-general}
  Let $F\in\D(H_{2,\alpha})$. Then $F^{(r)}\in\D(H_{2,2})$ for
  $r\geq2$.
  If $F\in\D(H_3)$, then $F^{(r)}\in\D(H_{2,2})$ for all $r\geq 1$.
\end{corollary}

\end{document}